\documentclass[a4paper,10pt,
colorlinks,urlcolor=black,linkcolor=black,citecolor=black,
]{article}
\usepackage{amsmath,amssymb}
\usepackage[dvips]{color}
\usepackage{cite}
\usepackage{hangcaption}
\usepackage{amsmath,amssymb}
\usepackage[dvips]{color}
\usepackage{cite}
\usepackage{eepic}
\usepackage{epic}
\usepackage{hangcaption}
\usepackage{amsmath}
\usepackage{amssymb}
\usepackage{amscd}
\usepackage{amsthm}
\usepackage{multicol}
\usepackage[dvipdfmx,%
 bookmarks=true,%
 bookmarksnumbered=true,%
 colorlinks=true,%
 setpagesize=false,%
 pdfkeywords={TeX; dvipdfmx; hyperref; color;}]{hyperref}
\def\cases{\left\{\begin{array}{ll}}
\def\endcases{\end{array}\right.}

\def\bigtimes{\mathop{\mbox{\Large $\times$}}}
\begin{document}
\setcounter{page}{1}
\vskip1.5cm
\begin{center}
{\Large \bf 
The confidence interval methods in quantum language
}
\vskip0.5cm
{\rm
\large
Shiro Ishikawa, Kohshi Kikuchi
}
\\
\vskip0.2cm
\rm
\it
Department of Mathematics, Faculty of Science and Technology,
Keio University,
\\ 
3-14-1, Hiyoshi, Kouhoku-ku Yokohama, Japan.
\\
E-mail:
ishikawa@math.keio.ac.jp,
kohshi.kikuchi@gmail.com
\end{center}
\par
\rm
\vskip0.3cm
\par
\noindent
{\bf Abstract}
\normalsize
\vskip0.5cm
\par
\noindent
Recently we proposed quantum language(or, measurement theory), which is characterized as the linguistic turn of the Copenhagen interpretation of quantum mechanics. Also, we consider that this is a kind of system theory such that it is applicable to both classical and quantum systems. As far as classical systems, it should be noted that quantum language is similar to statistics.
In this paper, we discuss the usual confidence interval methods in terms of quantum language. And we assert that three concepts
(i.e.,
"estimator" and "quantity" and "semi-distance)
are indispensable for the theoretical understanding of the
confidence interval methods.
Since our argument is quite elementary,
we hope that the readers acquire a new viewpoint of statistics, and agree that
our proposal is,
from the pure theoretical point of view, the true confidence interval methods.

\par
\noindent
(Key words: Confidence interval, Chi-squared distribution, Student's t-distribution)

\vskip1.0cm

\par

\def\Cal{\cal}
\def\bigstimes{\text{\large $\: \boxtimes \,$}}

\par
\noindent

\vskip0.2cm
\par
\noindent
\par
\noindent
\section{\large
Quantum language
(Axioms
and
Interpretation)
}
%

\rm
\par
\par
\noindent

In this section,
we shall mention the overview of quantum language
(or, measurement theory, in short, MT).
\par
\par
\rm
Quantum language is characterized as the linguistic turn of the Copenhagen interpretation of quantum mechanics({\it cf.} ref.
{{{}}}{\cite{Neum}}).
Quantum language (or, measurement theory ) has two simple rules
(i.e. Axiom 1(concerning measurement) and Axiom 2(concerning causal relation))
and the linguistic interpretation (= how to use the Axioms 1 and 2). 
That is,
\begin{align}
\underset{\mbox{(=MT(measurement theory))}}{\fbox{Quantum language}}
=
\underset{\mbox{(measurement)}}{\fbox{Axiom 1}}
+
\underset{\mbox{(causality)}}{\fbox{Axiom 2}}
+
\underset{\mbox{(how to use Axioms)}}{\fbox{linguistic interpretation}}
\label{eq1}
\end{align}
({\it cf.} refs.
{{{}}}{\cite{Ishi3}-\cite{Kikuchi}}).
\par
This theory is formulated in a certain $C^*$-algebra ${\cal A}$({\it cf.} ref.
{{{}}}{\cite{Saka}}), and is classified as follows:
\begin{itemize}
\item[(A)]
$
\quad
\underset{\text{\scriptsize }}{\text{MT}}
$
$\left\{\begin{array}{ll}
\text{quantum MT$\quad$(when ${\cal A}$ is non-commutative)}
\\
\\
\text{classical MT
$\quad$
(when ${\cal A}$ is commutative, i.e., ${\cal A}=C_0(\Omega)$)}
\end{array}\right.
$
\end{itemize}
where $C_0(\Omega)$
is
the $C^*$-algebra composed of all continuous 
complex-valued functions vanishing at infinity
on a locally compact Hausdorff space $\Omega$.

Since our concern in this paper is concentrated to 
the usual confidence interval
methods in statistics,
we devote ourselves to the commutative $C^*$-algebra $C_0(\Omega)$,
which is quite elementary.
Therefore, we believe that all statisticians
can understand our assertion
(i.e.,
a new viewpoint of the confidence interval
methods
).

Let $\Omega$ is a locally compact Hausdorff space, which is also called
a state space. And thus, an element $\omega (\in \Omega )$ is said to be a state.
Let $C(\Omega)$ be the $C^*$-algebra composed of all bounded continuous 
complex-valued functions on a locally compact Hausdorff space $\Omega$.
The norm $\| \cdot \|_{C(\Omega )}$ is usual, i.e.,
$\| f \|_{C(\Omega )} = \sup_{\omega \in \Omega } |f(\omega )|$
$(\forall f \in C(\Omega ))$.
\par
\rm
Motivated by Davies' idea ({\it cf.} ref.
{{{}}}{\cite{Davi}}) in quantum mechanics.
an observable ${\mathsf O}=(X, {\mathcal F}, F)$ in $C_0(\Omega )$
(or, precisely, in $C(\Omega )$) is defined as follows:
\begin{itemize}
\item[(B$_1$)]
$X$ is a topological space. 
${\mathcal F} ( \subseteq 2^X$(i.e., the power set of $X$) is a field,
that is, it satisfies the following conditions (i)--(iii):
(i):
$\emptyset \in {\cal F}$, 
(ii):$\Xi \in {\mathcal F} \Longrightarrow X\setminus \Xi  \in 
{\mathcal F}$,
(iii):
$\Xi_1, \Xi_2,\ldots, \Xi_n \in {\mathcal F} \Longrightarrow \cup_{k=1}^n \Xi_k \in {\mathcal F}$.
\item[(B$_2$)]
The map $F: {\cal F} \to C(\Omega )$ satisfies that 
\begin{align}
0 \le [F(\Xi )](\omega ) \le 1, \quad [F(X )](\omega )=1
\qquad
(\forall \omega \in \Omega )
\nonumber \end{align}
and moreover, 
if
\begin{align}
\Xi_1, \Xi_2,\ldots, \Xi_n, \ldots \in {\mathcal F},
\quad
\Xi_m \cap \Xi_n = \emptyset \quad( m \not= n ),
\quad
\Xi = \cup_{k=1}^\infty \Xi_k \in {\mathcal F},
\nonumber \end{align} 
then, it holds
\begin{align}
[F(\Xi)](\omega) = \lim_{n \to \infty } \sum_{k=1}^n [F(\Xi_k )](\omega )
\quad
(\forall \omega \in \Omega )
\nonumber \end{align}
\end{itemize}
Note that Hopf extension theorem
({\it cf.}
ref.
{{{}}}{\cite{Yosi}})
guarantees that
$(X, {\cal F}, [F(\cdot)](\omega))$
is regarded as the mathematical probability space.
\par
\noindent
\bf
Example 1
\rm
[Normal observable].
Let ${\mathbb R}$ be the set of the real numbers.
Consider the state space
$\Omega = {\mathbb R} \times {\mathbb R}_+$,
where 
${\mathbb R}_+=\{ \sigma \in {\mathbb R} | \sigma > 0 \}$.
Define the normal observable
${\mathsf O}_N = ({\mathbb R}, {\mathcal B}_{\mathbb R}, {{{N}}} )$ in $C_0({\mathbb R} \times {\mathbb R}_+)$ such that
\begin{align}
&
[{{{N}}}({\Xi})] ({} {}{\omega} {}) 
=
\frac{1}{{\sqrt{2 \pi }\sigma{}}}
\int_{{\Xi}} \exp[{}- \frac{({}{}{x} - {}{\mu}  {})^2 }{2 \sigma^2}    {}] d {}{x}
\label{eq2}
\\
&
\quad
({}\forall  {\Xi} \in {\cal B}_{{\mathbb R}{}}\mbox{(=Borel field in ${\mathbb R}$))},
\quad
\forall   {}{\omega} =(\mu, \sigma)   \in \Omega = {\mathbb R}{}\times {\mathbb R}_+).
\nonumber
\end{align}
In this paper, we devote ourselves to the normal observable.

\vskip0.5cm

%
%

\par
Now we shall briefly explain "quantum language (\ref{eq1})" in classical systems as follows:
A measurement of an observable
${\mathsf O}=(X, {\mathcal F}, F)$
for a system with a state $\omega (\in \Omega )$
is denoted by
${\mathsf M}_{C_0(\Omega)} ({\mathsf O}, S_{[\omega]})$.
By the measurement, a measured value $x (\in X)$ is obtained as follows:
\par
\noindent
\bf
Axiom 1
\rm
(Measurement)
\begin{itemize}
\item{}
\sl
The probability that a measured value $x$
$( \in X)$ obtained by the measurement 
${\mathsf{M}}_{{{C_0(\Omega)}}} ({\mathsf{O}}$
${ \equiv} (X, {\cal F}, F),$
{}{$ S_{[\omega_0]})$}
belongs to a set 
$\Xi (\in {\cal F})$ is given by
$
[F(\Xi) ](\omega_0 )
$.
\end{itemize}
\rm
\par
\noindent
\par
\noindent
\bf
Axiom 2
\rm
(Causality)
\begin{itemize}
\item{}
\sl
The causality is represented by a Markov operator
$\Phi_{21} : C_0(\Omega_2 ) \to C_0(\Omega_1 )$.
Particularly, the deterministic causality
is represented by a continuous map
$\pi_{21} : \Omega_1 \to \Omega_2$
\end{itemize}
\par
\noindent
\bf
Interpretation
\rm
(Linguistic interpretation).
Although there are several linguistic rules in quantum language, the following is the most important:
\begin{itemize}
\item{}
\sl
Only one measurement is permitted.
\end{itemize}
\rm
In order to read this paper,
it suffices to understand the above three.
\vskip0.5cm

Consider measurements ${\mathsf{M}}_{{{C_0(\Omega)}}} ({\mathsf{O}_k}$
${ \equiv} (X_k, {\cal F}_k, F_k),$
{}{$ S_{[\omega_0]})$}, $(k=1,2, \ldots, n )$. However, the linguistic interpretation says that
only one measurement is permitted.
Thus we must consider a simultaneous measurement or a parallel measurement.
The two are completely different, however in classical cases it suffices to consider only simultaneous measurement as follows.

\par
\noindent
\bf
Definition 1
\rm
[(i):Simultaneous observable].
\sl
Let
${\mathsf{O}_k}$
${ \equiv} (X_k, {\cal F}_k, F_k)$
$(k=1,2, \ldots, n)$
be an observable in $C_0(\Omega)$.
The simultaneous observable
$\bigtimes_{k=1}^{n}{\mathsf{O}_k}$
${ \equiv} (\bigtimes_{k=1}^n X_k, \bigstimes_{k=1}^n {\cal F}_k, \widehat{F}
(\equiv \bigtimes_{k=1}^n F_k
))$
in $C_0(\Omega)$ is defined by
\begin{align}
&
[\widehat{F}(\Xi_1 \times \cdots \times \Xi_n )](\omega)
(\equiv
[(\bigtimes_{k=1}^n F_k)(\Xi_1 \times \cdots \times \Xi_n )](\omega)
)
=
\bigtimes_{k=1}^n [F_k(\Xi_k)](\omega)
\label{eq3}
\\
&
\quad
\qquad
\qquad
(\forall \Xi_k \in {\mathcal F}_k \;\;(k=1, \ldots, n), \forall \omega \in \Omega )
\nonumber
\end{align}
Here,
$ \boxtimes_{k=1}^n {\cal F}_k$
is the smallest field including
the family
$\{
{\text{\large $\times$}}_{k=1}^n \Xi_k
$
$:$
$\Xi_k \in {\cal F}_k \; k=1,2,\ldots, n \}$.
If 
${\mathsf{O}}$
${ \equiv} (X, {\cal F}, F)$
is equal to
${\mathsf{O}_k}$
${ \equiv} (X_k, {\cal F}_k, F_k)$
$(k=1,2, \ldots, n)$,
then
the simultaneous observable
$\bigtimes_{k=1}^{n}{\mathsf{O}_k}$
${ \equiv} (\bigtimes_{k=1}^n X_k, \bigstimes_{k=1}^n {\cal F}_k, \widehat{F}
(\equiv \bigtimes_{k=1}^n F_k
))$
is denoted by
${\mathsf{O}^n}$
${ \equiv} (X^n, {\cal F}^n, F^n)$.
\par
\noindent
\rm
[(ii):Parallel observable].
\sl
Let
${\mathsf{O}_k}$
${ \equiv} (X_k, {\cal F}_k, F_k)$
be an observable in $C_0(\Omega_k)$,
$(k=1,2, \ldots, n)$.
The parallel observable
$\bigotimes_{k=1}^{n}{\mathsf{O}_k}$
${ \equiv} (\bigtimes_{k=1}^n X_k, \bigstimes_{k=1}^n {\cal F}_k, \widetilde{F}
(\equiv \bigotimes_{k=1}^n F_k
))$
in $C_0(\bigtimes_{k=1}^n \Omega_k)$ is defined by
\begin{align}
&
[\widetilde{F}(\Xi_1 \times \cdots \times \Xi_n )](\omega_1,\omega_2, \ldots,\omega_n)
(\equiv
[(\bigotimes_{k=1}^n F_k)(\Xi_1 \times \cdots \times \Xi_n )](\omega_1,\omega_2, \ldots,\omega_n)
)
=
\bigtimes_{k=1}^n [F_k(\Xi_k)](\omega_k)
\label{eq4}
\\
&
\quad
\qquad
\qquad
\quad
\qquad
\qquad
(\forall \Xi_k \in {\mathcal F}_k , \forall \omega_k \in \Omega_k,\;\;(k=1, \ldots, n))
\nonumber
\end{align}

\par
\noindent
\bf
Definition 2
\rm
[Image observable].
\sl
Let
${\mathsf{O}}$
${ \equiv} (X, {\cal F}, F)$
be observables in $C_0(\Omega)$.
The observable $f({\mathsf{O}})$ 
$({ \equiv} (Y, {\cal G}, G (\equiv F \circ f^{-1}))$
in $C_0(\Omega)$
is called the image observable of
${\mathsf{O}}$
by a map $f:X \to Y$,
if it holds that
\begin{align}
G( \Gamma ) = F( f^{-1}(\Gamma))
\qquad
(
\forall \Gamma \in {\mathcal G}
)
\label{eq5}
\end{align}

\par
\noindent
\bf
Example 2
\rm
[Simultaneous normal observable].
Let
${\mathsf O}_N = ({\mathbb R}, {\mathcal B}_{\mathbb R}, {{{N}}} )$ be the normal observable
in $C_0({\mathbb R} \times {\mathbb R}_+)$ in Example 1.
Let $n$ be a natural number.
Then, we get the simultaneous normal observable
${\mathsf O}_N^n = ({\mathbb R}^n, {\mathcal B}_{\mathbb R}^n, {{{N}}^n} )$ 
in $C_0({\mathbb R} \times {\mathbb R}_+)$.
That is,
\par
\noindent
\begin{align}
&
[{{{N}}}^n
(\bigtimes_{k=1}^n \Xi_k)]
({}\omega{})
=
\bigtimes_{k=1}^n
[{{{N}}}(\Xi_k)](\omega)
\nonumber
\\
=
&
\frac{1}{({{\sqrt{2 \pi }\sigma{}}})^n}
\underset{{\bigtimes_{k=1}^n \Xi_k }}{\int \cdots \int}
\exp[{}- \frac{\sum_{k=1}^n ({}{}{x_k} - {}{\mu}  {})^2 
}
{2 \sigma^2}    {}] d {}{x_1} d {}{x_2}\cdots dx_n
\label{eq6}
\\
&
\qquad 
({}\forall  \Xi_k \in {\cal B}_{{\mathbb R}{}}^{}
({}k=1,2,\ldots, n),
\quad
\forall   {}{\omega}=(\mu, \sigma )    \in \Omega = {\mathbb R}\times {\mathbb R}_+{}).
\nonumber
\end{align}

Consider the maps
$\overline{\mu}: {\mathbb R}^n \to {\mathbb R}$ 
and
${\overline{S}}: {\mathbb R}^n \to {\mathbb R}$ 
such that
\begin{align}
&
\overline{\mu}
(x) =
\overline{\mu}
(x_1,x_2,\ldots , x_n ) =
\frac{x_1 + x_2 + \cdots + x_n}{n}
\quad( \forall x=(x_1,x_2,\ldots , x_n ) \in {\mathbb R}^n )
\label{eq7}
\\
&
{{\overline{S}}}
(x) =
{{\overline{S}}}
(x_1,x_2,\ldots , x_n ) =
{\sum_{k=1}^n ( x_k - 
\overline{\mu}
(x))^2}
\quad( \forall x=(x_1,x_2,\ldots , x_n ) \in {\mathbb R}^n)
\label{eq8}
\end{align}
Thus, we have two
image observables
$\overline{\mu}({\mathsf O}_N^n) $
$= ({\mathbb R}, {\mathcal B}_{\mathbb R}, {{{N}}^n} \circ \overline{\mu}^{-1} )$
and
${{\overline{S}}}({\mathsf O}_N^n) $
$= ({\mathbb R}_+, {\mathcal B}_{{\mathbb R}_+}, {{{N}}^n} \circ {{\overline{S}}}^{-1} )$
in $C_0({\mathbb R} \times {\mathbb R}_+)$.

It is easy to see that
\begin{align}
&
[({{{N}}^n} \circ \overline{\mu}^{-1})(\Xi_1)](\omega)
=
\frac{1}{({{\sqrt{2 \pi }\sigma{}}})^n}
\underset{
\{ x \in {\mathbb R}^n \;:\; {\overline{\mu}}(x) \in \Xi_1 \}}
{\int \cdots \int}
\exp[{}- \frac{\sum_{k=1}^n ({}{}{x_k} - {}{\mu}  {})^2 
}
{2 \sigma^2}    {}] d {}{x_1} d {}{x_2}\cdots dx_n
\nonumber
\\
=
&
\frac{\sqrt{n}}{{\sqrt{2 \pi }\sigma{}}}
\int_{{\Xi_1}} \exp[{}- \frac{n({}{}{x} - {}{\mu}  {})^2 }{2 \sigma^2}    {}] d {}{x}
\label{eq9}
\intertext{and}
&
[({{{N}}^n} \circ {{{\overline{S}}}}^{-1})(\Xi_2)](\omega)
=
\frac{1}{({{\sqrt{2 \pi }\sigma{}}})^n}
\underset{
\{ x \in {\mathbb R}^n \;:\; {\overline{S}}(x) \in \Xi_2 \}}
{\int \cdots \int}
\exp[{}- \frac{\sum_{k=1}^n ({}{}{x_k} - {}{\mu}  {})^2 
}
{2 \sigma^2}    {}] d {}{x_1} d {}{x_2}\cdots dx_n
\nonumber
\\
=
&
\int_{\Xi_2 / \sigma^2} p^{{\chi}^2}_{n-1}({ x} ) {dx}
\label{eq10}
\\
&
\quad
({}\forall  {\Xi_1} \in {\cal B}_{{\mathbb R}{}},
\;\;
\forall \Xi_2 \in {\cal B}_{{\mathbb R}_+{}},
\quad
\forall   {}{\omega} =(\mu, \sigma)   \in \Omega \equiv {\mathbb R}{}\times {\mathbb R}_+).
\nonumber
\end{align}
Here, $p^{{\chi}^2}_{n-1}({ x} )$ is the chi-squared distribution with $n-1$ degrees of freedom. That is,
\begin{align}
p^{{\chi}^2}_{n-1}({ x} )
=
\frac{x^{(n-1)/2-1}e^{-x/2}}{2^{(n-1)/2} \Gamma ((n-1)/2)}
\quad ( x > 0)
\label{eq11}
\end{align}
where $\Gamma$ is the gamma function.

%
%

\section{
Fisher's maximum likelihood method}

\rm
\par
\noindent 
\par
It is usual to consider that
we do not know the pure state
$\omega_0$
$(
\in
\Omega
)$
when
we take a measurement
${\mathsf{M}}_{{{C_0(\Omega)}}} ({\mathsf{O}}, S_{[\omega_0]})$.
That is because
we usually take a measurement ${\mathsf{M}}_{{{C_0(\Omega)}}} ({\mathsf{O}},
S_{[\omega_0]})$
in order to know the state $\omega_0$.
Thus,
when we want to emphasize that
we do not know the state $\omega_0$,
${\mathsf{M}}_{{{C_0(\Omega)}}} ({\mathsf{O}}, S_{[\omega_0]})$
is denoted by
${\mathsf{M}}_{{{C_0(\Omega)}}} ({\mathsf{O}}, S_{[\ast]})$.
Also,
if
we know that a state $\omega_0$
belongs to
a certain set suitable $K$
$(\subseteq 
\Omega )$,
the
${\mathsf{M}}_{{{C_0(\Omega)}}} ({\mathsf{O}}, S_{[\omega_0]})$
is denoted by
${\mathsf{M}}_{{{C_0(\Omega)}}} ({\mathsf{O}}, S_{[\ast]}
(K) )$.

\vskip0.5cm
\par
\noindent
{\bf Theorem 1}
\rm
[Fisher's maximum likelihood method ({\it cf.} refs.
{{{}}}{\cite{Ishi4},\cite{Ishi5},\cite{Kikuchi}})].
\sl
Consider a measurement
${\mathsf M}_{{C_0(\Omega)}}
(
{\mathsf O}=(X , {\cal F} , F )
,$
$ S_{[*]}(K))$.
Assume that
we know that the measured value $x \;(\in X )$
obtained by a measurement
${\mathsf M}_{{C_0(\Omega)}}
(
{\mathsf O}=(X , {\cal F} , F )
,$
$ S_{[*]}(K))$
belongs to
$\Xi (\in {\cal F})$.
Then,
there is a reason to infer that the unknown state
$[\ast ]$ is equal to $\omega_0 (\in K )$ such that
\begin{align}
[F(\Xi)](\omega_0)
=
\max_{\omega \in K} [F(\Xi)](\omega)
\label{eq12}
\end{align}
if the righthand side of this formula exists.
Also, if $\Xi=\{x\}$, it suffices to calculate the $\omega_0 (\in K)$ 
such that
\begin{align}
\lim_{\Xi \supseteq \{x \}, \Xi \to \{x \}} \frac{[F(\Xi)](\omega_0)}{
\max_{\omega \in K} [F(\Xi)](\omega)}
=1
\label{eq13}
\end{align}
\par
\noindent
\rm
\bf
Example 3
\rm
[Fisher's maximum likelihood method].
Consider the simultaneous normal observable
${\mathsf O}_N^n = ({\mathbb R}^n, {\mathcal B}_{\mathbb R}^n, {{{N}}^n} )$
in $C_0({\mathbb R} \times {\mathbb R}_+)$ in the formula (\ref{eq6}).
Thus, we have the simultaneous measurement
${\mathsf M}_{C_0({\mathbb R} \times {\mathbb R}_+ )} ({\mathsf O}_N^n = ({\mathbb R}^n, {\mathcal B}_{\mathbb R}^n, {{{N}}^n} )$,
$S_{[\ast]}(K))$ 
in $C_0({\mathbb R} \times {\mathbb R}_+)$.
Assume that a measured value $x=(x_1, x_2, \ldots, x_n ) (\in
{\mathbb R}^n )$ is obtained by the measurement.
Since the likelihood function
$L_x(\mu, \sigma) $
is defined by
\begin{align}
&
\qquad 
L_x(\mu, \sigma)
=
\frac{1}{({{\sqrt{2 \pi }\sigma{}}})^n}
\exp[{}- \frac{\sum_{k=1}^n ({}{}{x_k} - {}{\mu}  {})^2 
}
{2 \sigma^2}    {}] 
\label{eq14}
\\
&
({}\forall x = (x_1, x_2, \ldots , x_n ) \in {\mathbb R}^n,
\quad
\forall   {}{\omega}=(\mu, \sigma )    \in \Omega = {\mathbb R}\times {\mathbb R}_+{}).
\nonumber
\end{align}
it suffices to calculate the following equations:
\begin{align}
\frac{\partial L_x(\mu, \sigma)}{\partial \mu}=0,
\quad
\frac{\partial L_x(\mu, \sigma)}{\partial \sigma}=0
\label{eq15}
\end{align}
Thus, Fisher's maximum likelihood method says as follows.
\par
\noindent
(i):
Assume that $K={\mathbb R} \times {\mathbb R}_+ $.
Solving the equation (\ref{eq15}), we can infer that
$[\ast]=(\mu, \sigma)$ 
$(\in 
{\mathbb R} \times {\mathbb R}_+
)$
such that
\begin{align}
\mu=\overline{\mu}(x) =\frac{x_1 + x_2+ \ldots + x_n }{n},
\quad
\sigma=
\overline{\sigma}(x)
=
\sqrt{\frac{{\overline{S}}(x)}{n}}=\sqrt{\frac{\sum_{k=1}^n (x_k - \overline{\mu}(x))^2}{n}}
\label{eq16}
\end{align}
\par
\noindent
(ii):
Assume that $K={\mathbb R} \times \{ \sigma_1 \}$
$( \subseteq {\mathbb R} \times {\mathbb R}_+)$.
It is easy to see that there is  a reason to infer that
$[\ast]=(\mu, \sigma)$ 
$(\in 
{\mathbb R} \times {\mathbb R}_+
)$
such that
\begin{align}
\mu=\overline{\mu}(x)=\frac{x_1 + x_2+ \ldots + x_n }{n},
\quad
\sigma= \sigma_1
\label{eq17}
\end{align}
\par
\noindent
(iii):
Assume that $K=\{\mu_1\} \times {\mathbb R}_+ $
$( \subseteq {\mathbb R} \times {\mathbb R}_+)$.
There is  a reason to consider that
$[\ast]=(\mu, \sigma)$ 
$(\in 
{\mathbb R} \times {\mathbb R}_+
)$
such that
\begin{align}
\mu=\mu_1
\quad
\sigma=
\sqrt{\frac{\sum_{k=1}^n (x_k - \mu_1)^2}{n}}
\label{eq18}
\end{align}

\vskip1.0cm
\par
\par
\rm
\section{Confidence interval }
\par
\noindent
\par

Let
${\mathsf O} = ({}X, {\cal F} , F{}){}$
be an observable
formulated in a commutative
$C^*$-algebra
${C_0(\Omega)}$.
Let $\Theta$ be a locally compact space with the 
semi-distance $d^x_{\Theta}$
$(\forall x \in X)$,
that is,
for each $x\in X$,
the map
$d^x_{\Theta}: \Theta^2 \to [0,\infty)$
satisfies that
(i):$d^x_\Theta (\theta, \theta )=0$,
(ii):$d^x_\Theta (\theta_1, \theta_2 )$
$=d^x_\Theta (\theta_2, \theta_1 )$,
(ii):$d^x_\Theta (\theta_1, \theta_3 )$
$\le d^x_\Theta (\theta_1, \theta_2 )
+
d^x_\Theta (\theta_2, \theta_3 )
$.
\noindent
\par
Let
$\pi:\Omega \to \Theta$
be a continuous map,
which is a kind of causal relation (in Axiom 2), and called
\it
{\lq\lq}quantity{\rq\rq}$.\; \;$
\rm
Let
$E:X \to \Theta$
be a continuous (or more generally, measurable) map,
which is called
\it
{\lq\lq}estimator{\rq\rq}$.\; \;$
\rm
Let
$\gamma$
be a real number such that
$0 \ll \gamma < 1$,
for example,
$\gamma = 0.95$.
For any state
$ \omega ({}\in \Omega)$,
define
the positive number
$\eta^\gamma_{\omega}$
$({}> 0)$
such that:
\begin{align}
\eta^\gamma_{\omega}
=
\inf
\{
\eta > 0:
[F(\{ x \in X \;:\; 
d^x_\Theta ( E(x) , \pi( \omega ) )
\le \eta
\}
)](\omega )
\ge \gamma
\}
\end{align}
For any
$x$
$({}\in X{})$,
put
\begin{align}
D_x^{\gamma}
=
\{
\pi({\omega})
(\in
\Theta)
:
\omega \in \Omega,
\;\;
d^x_\Theta ({}E(x),
\pi(\omega )
)
\le
\eta^\gamma_{\omega }
\}.
\label{eq20} 
\end{align}
The $D_x^{\gamma}$ is called
\it
the $({}\gamma{})$-confidence interval
of
$x$.
\rm
\par
Note that,
\begin{enumerate}
\item[(C)]
\it
for any
$\omega_0
({}\in
\Omega)$,
the probability,
that
the measured value $x$
obtained
by the measurement
${\mathsf M}_{C_0(\Omega)} \big({}{\mathsf O}:= ({}X, {\cal F} , F{})  ,$
$ S_{[\omega_0 {}] } \big)$
satisfies the following
condition $(\flat)$,
is larger than
$\gamma$
({}e.g., $\gamma= 0.95${}).
\begin{enumerate}
\item[$(\flat)$]
$\qquad$
$\qquad$
$  d^x_\Theta (E(x),  \pi(\omega_0){}) \le  {\eta }^\gamma_{\omega_0}  $.
\end{enumerate}
\end{enumerate}
\par
\noindent
Assume that
we get
a measured value
$x_0$
by
the measurement
${\mathsf M}_{C_0(\Omega)} \big({}{\mathsf O}:= ({}X, {\cal F} , F{})  ,$
$ S_{[\omega_0 {}] } \big)$.
Then,
we see the following equivalence:
\begin{align}
(\flat) \; \Longleftrightarrow \;
\;
D_{x_0}^\gamma
\ni
\pi( \omega_0).
\label{eq21} 
\end{align}

\par
\noindent
\begin{center}
\unitlength=0.4mm
\begin{picture}(230,75)
\put(-19,0){{
\put(40,16){\scriptsize $x_0$}
\qbezier(40,20)(100,61)(157,42)
\qbezier(153,32)(200,-5)(257,32)
\path(107,49)(115,48)(107,45)
\put(112,53){$E$}
\put(105,-35){
\put(102,55){$\pi$}
\path(107,51)(99,48)(107,45)}
\put(40,20){\circle*{1}}
\put(157,41){\circle*{1}}
\put(155,45){\scriptsize $E({}x_0)$}
\put(151,33){\scriptsize $\; \pi(\omega_0)$}
\put(251,30){\scriptsize $\;\;\; \cdot \; \omega_0$}
\put(149,34){ \circle*{1} }
\put(175,35){$ D_{x_0}^\gamma$}
\put(153,63){ $\Theta$}
\put(253,63){ $\Omega$}
\put(57,63){$X$}
\allinethickness{0.5mm}
\put(60,30){\oval(70,60)}
\put(160,30){\oval(70,60)}
\put(260,30){\oval(70,60)}
\allinethickness{0.3mm}
\put(157,42){\ellipse{30}{30}}
}}
\put(40,-20){\bf Figure 1.
\rm
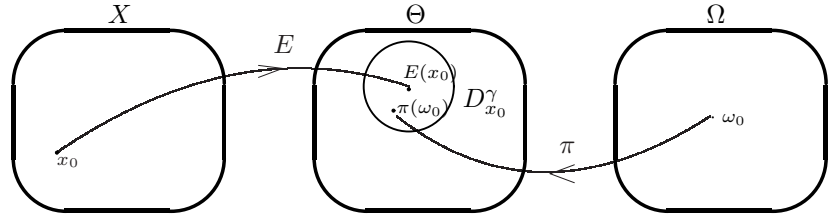
Confidence interval
$D^\gamma_{x_0}$
}
\end{picture}
\end{center}

\par
\vskip1.0cm
\par

%
%

\par
\vskip1.0cm
\par

\par
Summing the above argument,
we have the following proposition.
\par
\noindent
\bf 
\bf
Theorem 2
\rm
[{}Confidence interval{}].
\sl
Let
${\mathsf O} = ({}X, {\cal F} , F{}){}$
be an observable
formulated in a commutative
$C^*$-algebra
${C_0(\Omega)}$.
Let
$\omega_0$
be any fixed state,
i.e.,
$\omega_0 \in
\Omega
$,
Consider
a measurement
${\mathsf M}_{C_0(\Omega)} \big({}{\mathsf O}:= ({}X, {\cal F} , F{})  ,$
$ S_{[\omega_0 {}] } \big)$.
Let $\Theta$ be a locally compact space with the 
semi-distance $d^x_{\Theta}$
$(\forall x \in X)$.
Let
$\pi:\Omega \to \Theta$
be a quantity.
Let
$E:X \to \Theta$
be an estimator.
Let
$\gamma$
be a real number such that
$0 \ll \gamma < 1$,
for example,
$\gamma = 0.95$.
%
%
%
%
%
%
%
%
%
For any
$x
({}\in X{})$,
define
$D_x^{\gamma}$
as in (\ref{eq20}).
Then, we see,
\begin{enumerate}
\item[$(\sharp)$]
the probability
that
the measured value
$x_0
({}\in X)$
obtained
by the measurement
${\mathsf M}_{C_0(\Omega)} \big({}{\mathsf O}:= ({}X, {\cal F} , F{})  ,$
$ S_{[\omega_0 {}] } \big)$
satisfies
the condition
that
\begin{align}
\text{
$D_{x_0}^{\gamma} \ni \pi(\omega_0 )$
},
\label{eq22} 
\end{align}
is larger than
$\gamma$.
\end{enumerate}
\rm
This theorem is the generalization of our proposal in refs.{\cite{Ishi5} and \cite{Ishi10}}.
\par
\noindent
\bf
Remark 1
\rm
[The statistical meaning of Theorem 2].
Consider the simultaneous measurement
${\mathsf M}_{C_0(\Omega)} \big({}{\mathsf O}^J:= ({}X^J, {\cal F}^J , F^J{})  ,$
$ S_{[\omega_0 {}] } \big)$,
and assume that a measured value $x=(x_1,x_2, \ldots , x_J)( \in X^J)$ is obtained by the simultaneous measurement.
Then, it surely holds that
\begin{align}
\lim_{J \to \infty }
\frac{\mbox{Num} [\{ j \;|\; D_{x_j}^{\gamma} \ni \pi( \omega_0)]}{J}
\ge \gamma (= 0.95)
\label{eq23}
\end{align}
where
$\mbox{Num} [A]$ is the number of the elements of the set $A$.
Hence Theorem 2 can be tested by numerical analysis
(with random number).

\par
\section{Examples}
\subsection{The case that $\Omega=\Theta$,
and
$d^x_\Theta$ does not depend on $x$}
\par
\noindent
\rm
\par
In this section, we assume that  $\Omega=\Theta$,
that is, we do not need $\Theta$ but $\Omega$.
And moreover, we assume that
$d^x_\Theta$ does not depend on $x$.

The arguments in this section are continued from Example 2.
Consider the simultaneous measurement
${\mathsf M}_{C_0({\mathbb R} \times {\mathbb R}_+)}$
$({\mathsf O}_N^n = ({\mathbb R}^n, {\mathcal B}_{\mathbb R}^n, {{{N}}^n}) ,$
$S_{[(\mu, \sigma)]})$
in $C_0({\mathbb R} \times {\mathbb R}_+)$.
Thus,
we consider that
$\Omega = {\mathbb R} \times {\mathbb R}_+$,
$X={\mathbb R}^n$.
\rm
The formulas (\ref{eq7}) and (\ref{eq8}) urge us to
define the estimator 
$E: {\mathbb R}^n \to \Omega (\equiv \Theta  \equiv {\mathbb R} \times {\mathbb R}_+ )$
such that
\begin{align}
E(x)=E(x_1, x_2, \ldots , x_n )
=
(\overline{\mu}(x),
(\overline{\sigma}(x)
)
=
\Big(\frac{x_1 + x_2 + \cdots + x_n}{n},
\sqrt{\frac{\sum_{k=1}^n ( x_k - 
\overline{\mu}
(x))^2}{n}}
\Big)
\label{eq24}
\end{align}
Let
$\gamma$
be a real number such that
$0 \ll \gamma < 1$,
for example,
$\gamma = 0.95$.
\par
\noindent
\bf
Example 4
\rm
[Confidence interval for the semi-distance $d_{\Omega}^{(1)}$].
Consider the following semi-distance $d_{\Omega}^{(1)}$
in the state space ${\mathbb R} \times {\mathbb R}_+$:
\begin{align}
d_{\Omega}^{(1)}((\mu_1,\sigma_1), (\mu_2,\sigma_2))
=
|\mu_1 - \mu_2|
\label{eq25}
\end{align}

%

For any
$ \omega=(\mu, \sigma )  ({}\in\Omega=
{\mathbb  R} \times {\mathbb R}_+ )$,
define
the positive number
$\eta^\gamma_{\omega}$
$({}> 0)$
such that:
\begin{align}
\eta^\gamma_{\omega}
=
\inf
\{
\eta > 0:
[F ({}E^{-1} ({}
{{\rm Ball}_{d_\Omega^{(1)}}}(\omega ; \eta{}))](\omega )
\ge \gamma
\}
\nonumber
\end{align}
where
${{\rm Ball}_{d_\Omega^{(1)}}}(\omega ; \eta)$
$=$
$\{ \omega_1
({}\in\Omega):
d_\Omega^{(1)} ({}\omega, \omega_1{}) \le \eta \}$
$=
[\mu - \eta , \mu + \eta ]
\times {\mathbb R}_+$

Hence we see that
\begin{align}
&
E^{-1}({{\rm Ball}_{d_\Omega^{(1)}}}(\omega ; \eta ))
=
E^{-1}([\mu - \eta , \mu + \eta ] \times {\mathbb R}_+)
\nonumber
\\
=
&
\{
(x_1, \ldots , x_n )
\in {\mathbb R}^n
\;:
\;
\mu - \eta 
\le
\frac{x_1+\ldots + x_n }{n} \le  \mu + \eta 
\}
\label{eq26}
\end{align}
Thus,
\begin{align}
&
[{{{N}}}^n
(E^{-1}({{\rm Ball}_{d_\Omega^{(1)}}}(\omega ; \eta ))]
({}\omega{})
\nonumber
\\
=
&
\frac{1}{({{\sqrt{2 \pi }\sigma{}}})^n}
\underset{{
\mu - \eta 
\le
\frac{x_1+\ldots + x_n }{n} \le  \mu + \eta 
}}{\int \cdots \int}
\exp[{}- \frac{\sum_{k=1}^n ({}{}{x_k} - {}{\mu}  {})^2 
}
{2 \sigma^2}    {}] d {}{x_1} d {}{x_2}\cdots dx_n
\nonumber
\\
=
&
\frac{1}{({{\sqrt{2 \pi }\sigma{}}})^n}
\underset{{
 - \eta 
\le
\frac{x_1+\ldots + x_n }{n} \le   \eta 
}}{\int \cdots \int}
\exp[{}- \frac{\sum_{k=1}^n ({}{}{x_k}  {}{}  {})^2 
}
{2 \sigma^2}    {}] d {}{x_1} d {}{x_2}\cdots dx_n
\nonumber
\\
=
&
\frac{\sqrt{n}}{{\sqrt{2 \pi }\sigma{}}}
\int_{{- \eta}}^{\eta} \exp[{}- \frac{{n}{x}^2 }{2 \sigma^2}] d {x}
=
\frac{1}{{\sqrt{2 \pi }{}}}
\int_{{- \sqrt{n} \eta/\sigma}}^{\sqrt{n} \eta / \sigma} \exp[{}- \frac{{x}^2 }{2 }] d {x}
\label{eq27}
\end{align}
Solving the following equation:
\begin{align}
\frac{1}{{\sqrt{2 \pi }{}}}
\int^{-z((1-\gamma)/2)}_{-\infty} \exp[{}- \frac{{x}^2 }{2 }] d {x}
=
\frac{1}{{\sqrt{2 \pi }{}}}
\int_{z((1-\gamma)/2)}^{\infty} \exp[{}- \frac{{x}^2 }{2 }] d {x}
=
\frac{1- \gamma}{2}
\label{eq28}
\end{align}
we define that
\begin{align}
\eta^\gamma_{\omega} = 
 \frac{\sigma}{\sqrt{n}}
 z(\frac{1-\gamma}{2})
\label{eq29}
\end{align}

Therefore,
for any
$x$
$({}\in {\mathbb R}^n)$,
we get $D_x^{\gamma}$
(
the $({}\gamma{})$-confidence interval
of
$x$
)
as follows:
\begin{align}
D_x^{\gamma}
&
=
\{
{\omega}
(\in
\Omega)
:
d_\Omega ({}E(x),
\omega)
\le
\eta^\gamma_{\omega }
\}
\nonumber
\\
&
=
\{ (\mu, \sigma ) \in {\mathbb R} \times {\mathbb R}_+
\;:\;
| \mu - \overline{\mu}(x)|
=
| \mu - \frac{x_1+ \ldots + x_n}{n}|
\le 
 \frac{\sigma}{\sqrt{n}}
 z(\frac{1-\gamma}{2})
 \}
\label{eq30}
\end{align}
\rm

\par
\noindent
\begin{center}
\unitlength=0.43mm
\begin{picture}(250,150)
\put(75,50)
{{
\put(100,-10){${\mathbb R}$}
\put(-8,85){${\mathbb R}_+$}
\put(12,50){\large{${D_x^{\gamma}}$}}
\put(-30,0){\vector(1,0){130}}
\put(-15,0){\vector(0,1){90}}
\allinethickness{0.1mm}
\multiput(-43,58)(3,-3){20}{\line(0,1){30}}
\multiput(74,58)(-3,-3){20}{\line(0,1){30}}
\path(15,-2)(15,2)
\put(14,-8){$\overline{\mu}(x)$}
\put(25,50){}
\thicklines
\put(-50,65){\line(1,-1){65}}
\put(81,65){\line(-1,-1){65}}
}}
\put(10,20){\bf Figure 2.
\rm
Confidence interval
$D_x^{\gamma}$
for the semi-distance $d_{\Omega}^{(1)}$
%
}
\end{picture}
\end{center}
Thus, strictly speaking, the "confidence interval" should be said to be
the "confidence domain" in quantum language.


\par
\noindent
\bf
Example 5
\rm
[Confidence interval for the semi-distance $d_{\Omega}^{(2)}$].
\rm
Consider the following semi-distance $d_{\Omega}^{(2)}$
in ${\mathbb R} \times {\mathbb R}_+$:
\begin{align}
d_{\Omega}^{(2)}((\mu_1,\sigma_1), (\mu_2,\sigma_2))
=
|
\int_{\sigma_1}^{\sigma_2} \frac{1}{\sigma} d \sigma
|
=
|\log{\sigma_1} - \log{\sigma_2}|
\label{eq31}
\end{align}

For any
$ \omega=(\mu, \sigma )  ({}\in\Omega=
{\mathbb  R} \times {\mathbb R}_+ )$,
define
the positive number
$\eta^\gamma_{\omega}$
$({}> 0)$
such that:
\begin{align}
\eta^\gamma_{\omega}
=
\inf
\{
\eta > 0:
[F ({}E^{-1} ({}
{{\rm Ball}_{d_\Omega^{(2)}}}(\omega ; \eta{}))](\omega )
\ge \gamma
\}
\label{eq32}
\end{align}
where
${{\rm Ball}_{d_\Omega^{(2)}}}(\omega ; \eta)$
$=$
$\{ \omega_1
({}\in\Omega):
d_\Omega^{(2)} ({}\omega, \omega_1{}) \le \eta \}$.
Note that
\begin{align}
{{\rm Ball}_{d_\Omega^{(2)}}}(\omega ; \eta )
=
{{\rm Ball}_{d_\Omega^{(2)}}}((\mu ; \sigma ), \eta )
=
{\mathbb R} \times \{ \sigma' \;:\; |\log(\sigma'/\sigma)| \le \eta
\}
=
{\mathbb R} \times [\sigma e^{-\eta} , \sigma e^{\eta} ]
\label{eq33}
\end{align}
Then,
\begin{align}
&
E^{-1}( {{\rm Ball}_{d_\Omega^{(2)}}}(\omega ; \eta ))
=
E^{-1}({\mathbb R} \times 
[\sigma e^{-\eta} , \sigma e^{\eta} ]
)
\nonumber
\\
=
&
\{
(x_1, \ldots , x_n )
\in
{\mathbb R}^n
\;:
\;
\sigma e^{- \eta }
\le
\Big(
\frac{\sum_{k=1}^n ( x_k - 
\overline{\mu}
(x))^2}{n}
\Big)^{1/2}
\le
\sigma e^{ \eta }
\}
\label{eq34}
\end{align}
Hence we see, by (\ref{eq10}), that
\begin{align}
&
[{{{N}}}^n
(E^{-1}({{\rm Ball}_{d_\Omega^{(2)}}}(\omega; \eta ))]
({}\omega{})
\nonumber
\\
=
&
\frac{1}{({{\sqrt{2 \pi }\sigma{}}})^n}
\underset{{
\sigma^2 e^{- 2 \eta }
\le
\frac{\sum_{k=1}^n ( x_k - 
\overline{\mu}
(x))^2}{n}
\le
\sigma^2 e^{ 2 \eta }
}}{\int \cdots \int}
\exp[{}- \frac{\sum_{k=1}^n ({}{}{x_k} - {}{\mu}  {})^2 
}
{2 \sigma^2}    {}] d {}{x_1} d {}{x_2}\cdots dx_n
\nonumber
\\
=
&
\int_{{n}  e^{- 2 \eta}}^{{n}  e^{ 2 \eta}} 
p^{\chi^2}_{n-1} (x ) 
dx
\label{eq35}
\end{align}
Using the chi-squared distribution $p^{{\chi}^2}_{n-1}({ x} )$
(with $n-1$ degrees of freedom) in (\ref{eq11}),
define the $\eta^\gamma_{\omega}$ such that
\begin{align}
\gamma
=
\int_{{n} e^{-2 \eta^\gamma_{\omega}}}^{{n}  e^{2 \eta^\gamma_{\omega}}} 
p^{\chi^2}_{n-1} (x ) dx
\label{eq36}
\end{align}
where it should be noted that
the $\eta^\gamma_{\omega}$
depends on only $\gamma$ and $n$.
Thus, put
\begin{align}
\eta^\gamma_{\omega} = \eta^\gamma_{n}
\label{eq37}
\end{align}
Hence we get,
for any
$x$
$({}\in X{})$,
the
$D_x^{\gamma}$
(
the $({}\gamma{})$-confidence interval
of
$x$
)
as follows:
\begin{align}
D_x^{\gamma}
&
=
\{
{\omega}
(\in
\Omega)
:
d^{(2)}_\Omega ({}E(x),
\omega)
\le
\eta^\gamma_{n}
\}
\nonumber
\\
&
=
\{ (\mu, \sigma ) \in {\mathbb R} \times {\mathbb R}_+
\;:
\;
\sigma e^{- \eta^\gamma_{n} }
\le
\Big(
\frac{\sum_{k=1}^n ( x_k - 
\overline{\mu}
(x))^2}{n}
\Big)^{1/2}
\le
\sigma e^{ \eta^\gamma_{n} }
\}
\label{eq38}
\intertext{Recalling (\ref{eq16}), i.e.,
$\overline{\sigma}(x)
=
\Big(
\frac{\sum_{k=1}^n ( x_k - 
\overline{\mu}
(x))^2}{n}
\Big)^{1/2}
={(\frac{{\overline{S}}(x)}{n})}^{1/2}
$,
we conclude that}
D_x^{\gamma}
&
=
\{ (\mu, \sigma ) \in {\mathbb R} \times {\mathbb R}_+
\;:
\;
\overline{\sigma}(x)
e^{ - \eta^\gamma_{n}}
\le
\sigma
\le
\overline{\sigma}(x)
e^{ \eta^\gamma_{n}}
\}
\nonumber
\\
&
=
\{ (\mu, \sigma ) \in {\mathbb R} \times {\mathbb R}_+
\;:
\;
\frac{e^{ - 2\eta^\gamma_{n}}}{n}{\overline{S}}(x)
\le
\sigma^2
\le
\frac{e^{ 2\eta^\gamma_{n}}}{n}{\overline{S}}(x)
\}
\label{eq39}
\end{align}

%

\par
\noindent
\begin{center}
\unitlength=0.43mm
\begin{picture}(250,150)
\put(75,50)
{{
\put(100,-10){${\mathbb R}$}
\put(-8,86){${\mathbb R}_+$}
\put(90,36){\large{${D_x^{\gamma}}$}}
\put(-30,0){\vector(1,0){130}}
\put(-15,0){\vector(0,1){90}}
\put(25,50){}
\thicklines
\allinethickness{0.1mm}
\multiput(-30,25)(4,0){30}{\line(0,1){30}}
\put(-33,56){\line(1,0){135}}
\put(-33,25){\line(1,0){135}}
\put(-33,65){
\put(30,5){\vector(-1,-1){10}}
\put(30,5){
${\overline{\sigma} (x)
e^{\eta^\gamma_{n} }
}
$
}
}
\put(-33,5){
\put(30,5){\vector(-1,1){10}}
\put(30,5){
${\overline{\sigma} (x)
e^{- \eta^\gamma_{n}}
}
$
}
}
}}
\put(20,25){\bf Figure 3.
\rm
Confidence interval
$D_x^{\gamma}$
for the semi-distance $d_{\Omega}^{(2)}$
}
\end{picture}
\end{center}
%
For example, in the case that $n=3$, $\gamma=0.95$, the (\ref{eq36}) says that
\begin{align}
0.95=\gamma
&
=
\int_{{3} e^{-2 \eta^\gamma_{n}}}^{{3}  e^{2 \eta^\gamma_{n}}} 
p^{\chi^2}_{2} (x ) dx
=
\int_{{3} e^{-2 \eta^\gamma_{n}}}^{{3}  e^{2 \eta^\gamma_{n}}} 
\frac{e^{-x/2}}{2^{2/2} \Gamma (1)}
dx
=
\Big[{-e^{-x/2}}\Big]_{x={{3}  e^{-2 \eta^\gamma_{n}}} }^{x={{3}  e^{2 \eta^\gamma_{n}}} }
\nonumber
\\
&
=
e^{-{\frac{3}{2}  e^{-{2} \eta^\gamma_{n}}}}-e^{-{\frac{3}{2}  e^{2 \eta^\gamma_{n}}}}
\label{eq40}
\end{align}
which implies that
%
\begin{align}
\quad
e^{- \eta^{0.95}_{3}}=0.1849,
\qquad
e^{ \eta^{0.95}_{3}}
=5.4077,
\label{eq41}
\end{align}
and,
\begin{align}
\quad
e^{ -2 \eta^{0.95}_{3}}/3
=0.0114\cdots,
\quad
e^{ 2 \eta^{0.95}_{3}}/3
=9.748 \cdots
\label{eq42}
\end{align}
Thus, we see that
\begin{align}
D_x^{0.95}
=
\{ (\mu, \sigma ) \in {\mathbb R} \times {\mathbb R}_+
\;:
\;
(0.0114\cdots)\cdot {\overline{S}}(x)
\le
\sigma^2
\le
(9.748\cdots)\cdot {\overline{S}}(x)
\}
\label{eq43}
\end{align}
\par
\noindent
\bf
Remark 2.
\rm
[Other estimator].
Instead of (\ref{eq24}), we may consider 
the unbiased estimator 
$E': {\mathbb R}^n \to \Omega ( \equiv {\mathbb R} \times {\mathbb R}_+ )$
such that
\begin{align}
E'(x)=E(x_1, x_2, \ldots , x_n )
=
(\overline{\mu}(x),
(\overline{\sigma}'(x)
)
=
\Big(\frac{x_1 + x_2 + \cdots + x_n}{n},
\sqrt{\frac{\sum_{k=1}^n ( x_k - 
\overline{\mu}
(x))^2}{n-1}}
\Big)
\label{eq44}
\end{align}
In this case,
we see that
\begin{align}
(D_x^{\gamma})'
&
=
\{ (\mu, \sigma ) \in {\mathbb R} \times {\mathbb R}_+
\;:
\;
\overline{\sigma}'(x)
e^{ - (\eta^\gamma_{n})'}
\le
\sigma
\le
\overline{\sigma}'(x)
e^{ (\eta^\gamma_{n})'}
\}
\nonumber
\\
&
=
\{ (\mu, \sigma ) \in {\mathbb R} \times {\mathbb R}_+
\;:
\;
\frac{e^{ - 2(\eta^\gamma_{n})'}}{n-1}{\overline{S}}(x)
\le
\sigma^2
\le
\frac{e^{ 2(\eta^\gamma_{n})'}}{n-1}{\overline{S}}(x)
\}
\label{eq45}
\end{align}
where
the $(\eta^\gamma_{n})'$
is defined by
\begin{align}
\gamma
=
\int_{{(n-1)} e^{-2(\eta^\gamma_{n})'}}^{{(n-1)}  e^{2(\eta^\gamma_{n})'}} 
p^{\chi^2}_{n-1} (x ) dx
\label{eq46}
\end{align}

For example, in the case that $n=3$, $\gamma=0.95$, the (\ref{eq36}) says that
\begin{align}
0.95=\gamma
&
=
\int_{{2} e^{-2 (\eta^\gamma_{n})'}}^{{}  2e^{2 (\eta^\gamma_{n})'}} 
p^{\chi^2}_{2} (x ) dx
=
\int_{{2} e^{-2 (\eta^\gamma_{n})'}}^{{2}  e^{2 (\eta^\gamma_{n})'}} 
\frac{e^{-x/2}}{2^{2/2} \Gamma (1)}
dx
=
\Big[{-e^{-x/2}}\Big]_{x={{2}  e^{-2 (\eta^\gamma_{n})'}} }^{x={{2}  e^{2 (\eta^\gamma_{n})'}} }
\nonumber
\\
&
=
e^{-{  e^{-{2} (\eta^\gamma_{n})'}}}-e^{-{  e^{2 (\eta^\gamma_{n})'}}}
\label{eq47}
\end{align}
which implies that
%
\begin{align}
\quad
e^{- (\eta^{0.95}_{3})'}=0.2265,
\quad
e^{ (\eta^{0.95}_{3})'}
=4.4154
\label{eq48}
\end{align}
Thus,
\begin{align}
e^{ -2 (\eta^{0.95}_{3})'}/2
=0.00256\cdots
\qquad
e^{ 2 (\eta^{0.95}_{3})'}/2
=9.748 \cdots
\quad
\label{eq49}
\end{align}
Thus, we see that
\begin{align}
(D_x^{0.95})'
=
\{ (\mu, \sigma ) \in {\mathbb R} \times {\mathbb R}_+
\;:
\;
(0.00256\cdots)\cdot {\overline{S}}(x)
\le
\sigma^2
\le
(9.748\cdots)\cdot {\overline{S}}(x)
\}
\label{eq50}
\end{align}
Hence it should be noted that
$
D_x^{\gamma}\not= (D_x^{\gamma})'
$.
\par
\vskip0.3cm
\par
\noindent
\bf
Remark 3
\rm
[Other semi-distance $d_{\Omega}^{(3)}$].
We believe that
the semi-distance $d_{\Omega}^{(2)}$
is natural in Example 5, although
we have no firm reason to believe in it.
For example, consider a positive continuous function
$h: {\mathbb R}_+ \to {\mathbb R}_+$.
%
Then,
we can define
another semi-distance $d_{\Omega}^{(3)}$
in the state space ${\mathbb R} \times {\mathbb R}_+$:
\begin{align}
d_{\Omega}^{(3)}((\mu_1,\sigma_1), (\mu_2,\sigma_2))
=
| \int_{\sigma_1}^{\sigma_2} h(\sigma) d \sigma
|
\label{eq51}
\end{align}
Thus, many ($\gamma$)-confidence intervals exist,
though the $\eta_n^\gamma$ may depend on $\omega$.
Now, we have the following problem:
\begin{itemize}
\item{}
Is there a better $h(\sigma )$ than the $1/ \sigma $?
\end{itemize}
whose answer we do not know.
\par
\vskip0.5cm
\noindent
\bf
Remark 4
\rm
[So called $\alpha$-point method].
In many books, it conventionally is recommended as follows:
\begin{align}
(D_x^{\gamma})''
&
=
\Big\{ (\mu, \sigma ) \in {\mathbb R} \times {\mathbb R}_+
\;:
\;
\sqrt{
\frac{\sum_{k=1}^n (x_k - \overline{\mu}(x) )^2}{\chi_\infty^2}
}
\le
\sigma
\le
\sqrt{
\frac{\sum_{k=1}^n (x_k - \overline{\mu}(x) )^2}{\chi_0^2}
}
\Big\}
\nonumber
\\
&
=
\Big\{ (\mu, \sigma ) \in {\mathbb R} \times {\mathbb R}_+
\;:
\;
\frac{{\overline{S}}(x)}{{\chi_\infty^2}}
\le
\sigma^2
\le
\frac{{\overline{S}}(x)}{{\chi_0^2}}
%
\Big\}
\label{eq52}
\end{align}
where
\begin{align}
&
\int_{0}^{\chi_0^2}p_{n-1}^{\chi^2}(x) dx =
\int_{\chi_\infty^2}^\infty p_{n-1}^{\chi^2}(x) dx
=(1- \gamma)/2
\label{eq53}
\end{align}
which may be an analogy of (19).
\par
\noindent
\par
In the case that $n=3$, $\gamma =0.95$, we see
\begin{align}
&
\int_{0}^{0.0506}p_{2}^{\chi^2}(x) dx =
\int_{7.378}^\infty p_{2}^{\chi^2}(x) dx
=0.025
\label{eq54}
\intertext{and thus,}
&
{\frac{1}{\chi_\infty^2}}
=
{\frac{1}{7.378}}
=0.1355,
\quad
{\frac{1}{\chi_0^2}}
=
{\frac{1}{0.0506}}=19.763, \qquad 
\label{eq55}
\end{align}
Thus, we see that
\begin{align}
(D_x^{0.95})''
=
\{ (\mu, \sigma ) \in {\mathbb R} \times {\mathbb R}_+
\;:
\;
(0.1355\cdots) \cdot {\overline{S}}(x)
\le
\sigma^2
\le
(19.763\cdots) \cdot {\overline{S}}(x)
\}
\label{eq56}
\end{align}
which should be compared to the estimations (\ref{eq43}) and (\ref{eq50}).
It should be noted that
both estimator and semi-distance are not declared in 
this $\alpha$-point method.
Thus, we have the following problem:
\begin{itemize}
\item[(C)]
What is the $\alpha$-point method (\ref{eq52})?
\end{itemize}
This will be answered in the following remark.

\rm
\noindent
\vskip0.5cm
\par
\noindent
\bf
Remark 5
\rm
[
What is the $\alpha$-point method (\ref{eq52})?
].
Instead of (\ref{eq24}) or (\ref{eq44}), we consider 
the estimator 
$E'': {\mathbb R}^n \to \Omega ( \equiv {\mathbb R} \times {\mathbb R}_+ )$
such that
\begin{align}
E''(x)=E(x_1, x_2, \ldots , x_n )
=
(\overline{\mu}(x),
(\overline{\sigma}''(x)
)
=
\Big(\frac{x_1 + x_2 + \cdots + x_n}{n},
\sqrt{\frac{\sum_{k=1}^n ( x_k - 
\overline{\mu}
(x))^2}{cn}}
\Big)
\label{eq144}
\end{align}
where $c > 0$.
In this case,
by the same argument of (\ref{eq35}),
we see that
Then,
\begin{align}
&
(E'')^{-1}( {{\rm Ball}_{d_\Omega^{(2)}}}(\omega ; \eta ))
=
E^{-1}({\mathbb R} \times 
[\sigma e^{-\eta} , \sigma e^{\eta} ]
)
\nonumber
\\
=
&
\{
(x_1, \ldots , x_n )
\in
{\mathbb R}^n
\;:
\;
\sigma e^{- \eta }
\le
\Big(
\frac{\sum_{k=1}^n ( x_k - 
\overline{\mu}
(x))^2}{cn}
\Big)^{1/2}
\le
\sigma e^{ \eta }
\}
\label{eq134}
\end{align}
Hence we see, by (\ref{eq10}), that
\begin{align}
&
[{{{N}}}^n
(E^{-1}({{\rm Ball}_{d_\Omega^{(2)}}}(\omega; \eta ))]
({}\omega{})
\nonumber
\\
=
&
\frac{1}{({{\sqrt{2 \pi }\sigma{}}})^n}
\underset{{
\sigma^2 e^{- 2 \eta }
\le
\frac{\sum_{k=1}^n ( x_k - 
\overline{\mu}
(x))^2}{cn}
\le
\sigma^2 e^{ 2 \eta }
}}{\int \cdots \int}
\exp[{}- \frac{\sum_{k=1}^n ({}{}{x_k} - {}{\mu}  {})^2 
}
{2 \sigma^2}    {}] d {}{x_1} d {}{x_2}\cdots dx_n
\nonumber
\\
=
&
\int_{{cn}  e^{- 2 \eta}}^{{cn}  e^{ 2 \eta}} 
p^{\chi^2}_{n-1} (x ) 
dx
\label{eq235}
\end{align}
Hence we get,
for any
$x$
$({}\in X{})$,
the
$D_x^{\gamma}$
(
the $({}\gamma{})$-confidence interval
of
$x$
)
as follows:
\begin{align}
D_x^{\gamma}
&
=
\{
{\omega}
(\in
\Omega)
:
d^{(2)}_\Omega ({}E(x),
\omega)
\le
\eta^\gamma_{n}
\}
\nonumber
\\
&
=
\{ (\mu, \sigma ) \in {\mathbb R} \times {\mathbb R}_+
\;:
\;
\sigma e^{- \eta^\gamma_{n} }
\le
\Big(
\frac{\sum_{k=1}^n ( x_k - 
\overline{\mu}
(x))^2}{cn}
\Big)^{1/2}
\le
\sigma e^{ \eta^\gamma_{n} }
\}
\nonumber
\\
&
=
\{ (\mu, \sigma ) \in {\mathbb R} \times {\mathbb R}_+
\;:
\;
\;
cn \sigma^2 e^{- 2\eta^\gamma_{n} }
\le
\overline{S}(x)
\le
cn \sigma^2 e^{ 2\eta^\gamma_{n} }
%
\}
\nonumber
\\
&
=
\{ (\mu, \sigma ) \in {\mathbb R} \times {\mathbb R}_+
\;:
\;
\;
\frac{\overline{S}(x)}{
cn  e^{ 2\eta^\gamma_{n} }}
\le
\sigma^2
\le
\frac{\overline{S}(x)}{
cn  e^{- 2\eta^\gamma_{n} }}
%
%
\}
%
\label{eq139}
\end{align}

%
%
\par
\noindent
Using $\chi_0^2$
and
$\chi_\infty^2$
defined in (\ref{eq53}),
we obtain the following equation:
\begin{align*}
{cn}  e^{- 2 \eta^\gamma_{n}}= \chi_0^2,
\qquad
{cn}  e^{ 2 \eta^\gamma_{n}}= \chi_\infty^2
\end{align*}
Thus, it suffices to put
\begin{align}
c=\frac{\sqrt{\chi_0^2\cdot \chi_\infty^2}}{n}
\label{100}
\end{align}
in the estimator $E''$ of
(\ref{eq144}).
In this sense,
the $\alpha$-point method (\ref{eq52})
is true
({\it cf.}
Remark 1),
though it may be unnatural.

\subsection{The case that $\pi(\mu_1, \mu_2)=\mu_1- \mu_2$,
and
$d^x_\Theta$ does not depend on $x$
}
\par
\noindent
\rm
\par
The arguments in this section are continued from Example 2.
\par
\noindent
\bf
Example 6
\rm
\rm
[Confidence interval the the case that "$\pi(\mu_1, \mu_2)=\mu_1- \mu_2$"].
\rm
Consider the parallel measurement
${\mathsf M}_{C_0(({\mathbb R} \times {\mathbb R}_+) \times ({\mathbb R} \times {\mathbb R}_+))}$
$({\mathsf O}_N^n \otimes {\mathsf O}_N^m= ({\mathbb R}^n \times {\mathbb R}^m \ , {\mathcal B}_{\mathbb R}^n \bigstimes {\mathcal B}_{\mathbb R}^m, {{{N}}^n}
\otimes  {{{N}}^m}) ,$
$S_{[(\mu_1, \sigma_1, \mu_2 , \sigma_2)]})$
in $C_0(({\mathbb R} \times {\mathbb R}_+) \times ({\mathbb R} \times {\mathbb R}_+))$.
\rm

Assume that $\sigma_1$ and $\sigma_2$ are fixed and known. Thus, this parallel 
measurement is represented by
${\mathsf M}_{C_0({\mathbb R} \times{\mathbb R} )}$
$({\mathsf O}_{N_{\sigma_1}}^n \otimes {\mathsf O}_{N_{\sigma_1}}^m= ({\mathbb R}^n \times {\mathbb R}^m \ , {\mathcal B}_{\mathbb R}^n \bigstimes {\mathcal B}_{\mathbb R}^m, {{{N_{\sigma_1}}}^n}
\otimes  {{{N_{\sigma_2}}}^m}) ,$
$S_{[(\mu_1,  \mu_2 )]})$
in $C_0({\mathbb R} \times {\mathbb R} )$. Here, recall the (\ref{eq2}),
i.e.,
\par
\noindent
\begin{align}
&
[{{{N_\sigma}}}({\Xi})] ({} {}{\mu} {}) 
=
\frac{1}{{\sqrt{2 \pi }\sigma{}}}
\int_{{\Xi}} \exp[{}- \frac{({}{}{x} - {}{\mu}  {})^2 }{2 \sigma^2}    {}] d {}{x}
\quad
({}\forall  {\Xi} \in {\cal B}_{{\mathbb R}{}}\mbox{(=Borel field in ${\mathbb R}$))},
\quad
\forall  \mu   \in {\mathbb R}).
\label{eq57}
\end{align}
Therefore, we have the state space
$\Omega ={\mathbb R}^2
=
\{
\omega=(\mu_1, \mu_2) \;:\; \mu_1,\mu_2
 \in {\mathbb R} 
\}$.
Put $\Theta={\mathbb R}$ with the distance $d_\Theta ( \theta_1, \theta_2 )= |\theta_1-\theta_2|$
and consider the quantity $\pi:{\mathbb R}^2 \to
{\mathbb R}$ by
\begin{align}
\pi (\mu_1, \mu_2)= \mu_1-\mu_2
\label{eq58}
\end{align}
The estimator $E: \widehat{X}(=X \times Y =
{{\mathbb R}^n \times {\mathbb R}^m}) 
\to \Theta(={\mathbb R})$ is defined by
\begin{align}
E(x_1, \ldots, x_n,y_1, \ldots, y_m)
=
\frac{\sum_{k=1}^n x_k}{n}
-
\frac{\sum_{k=1}^m y_k}{m}
\label{eq59}
\end{align}

For any
$ \omega=(\mu_1, \mu_2 )  ({}\in \Omega=
{\mathbb  R} \times {\mathbb R} )$,
define
the positive number
$\eta^\gamma_{\omega}$
$({}> 0)$
such that:
\begin{align}
\eta^\gamma_{\omega}
=
\inf
\{
\eta > 0:
[F ({}E^{-1} ({}
{{\rm Ball}_{d_\Theta}}(\pi(\omega) ; \eta{}))](\omega )
\ge \gamma
\}
\nonumber
\end{align}
where
${{\rm Ball}_{d_\Theta}}(\pi(\omega) ; \eta)$
$=
[\mu_1 - \mu_2  - \eta , \mu_1 - \mu_2 + \eta ]$

Now let us calculate the $\eta^\gamma_{\omega}$ as follows:
\begin{align}
&
E^{-1}({{\rm Ball}_{d_\Theta}}(\pi(\omega) ; \eta ))
=
E^{-1}([\mu_1 - \mu_2 - \eta , \mu_1 - \mu_2 + \eta ] )
\nonumber
\\
=
&
\{
(x_1, \ldots , x_n, y_1, \ldots, y_m )
\in {\mathbb R}^n \times {\mathbb R}^m
\;:
\;
\mu_1 - \mu_2 - \eta 
\le
\frac{\sum_{k=1}^n x_k}{n}
-
\frac{\sum_{k=1}^m y_k}{m}
\le  \mu_1 - \mu_2 + \eta 
\}
\nonumber
\\
=
&
\{
(x_1, \ldots , x_n, y_1, \ldots, y_m )
\in {\mathbb R}^n \times {\mathbb R}^m
\;:
\;
- \eta 
\le
\frac{\sum_{k=1}^n (x_k - \mu_1)}{n}
-
\frac{\sum_{k=1}^m (y_k- \mu_2)}{m}
\le  \eta 
\}
\label{eq60}
\end{align}

Thus,
\begin{align}
&
[
({{{N_{\sigma_1}}}}^n
\otimes
{{{N_{\sigma_2}}}}^m
)
(E^{-1}({{\rm Ball}_{d_\Theta}}(\pi(\omega) ; \eta ))]
({}\omega{})
\label{eq61}
\\
=
&
\frac{1}{({{\sqrt{2 \pi }\sigma_1{}}})^n({{\sqrt{2 \pi }\sigma_2{}}})^m}
\nonumber
\\
&
\bigtimes 
\!\!\!\!\!\!
\underset{{
- \eta 
\le
\frac{\sum_{k=1}^n( x_k - \mu_1)}{n}
-
\frac{\sum_{k=1}^m (y_k- \mu_2)}{m}
\le  \eta 
%
%
%
}}{\int \cdots \int}
\exp[
{}- \frac{\sum_{k=1}^n ({}{}{x_k} - {}{\mu_1}  {})^2 
}
{2 \sigma_1^2}
{}- \frac{\sum_{k=1}^m ({}{}{y_k} - {}{\mu_2}  {})^2 
}
{2 \sigma_2^2}
] d {}{x_1} d {}{x_2}\cdots dx_nd {}{y_1} d {}{y_2}\cdots dy_m
\nonumber
\\
=
&
\frac{1}{({{\sqrt{2 \pi }\sigma_1{}}})^n({{\sqrt{2 \pi }\sigma_2{}}})^m}
\underset{{
- \eta 
\le
\frac{\sum_{k=1}^n x_k }{n}
-
\frac{\sum_{k=1}^m y_k}{m}
\le  \eta 
%
%
%
}}{\int \cdots \int}
\exp[
- \frac{
\sum_{k=1}^n {x_k}^2 
}
{2 \sigma_1^2}
- \frac{
\sum_{k=1}^m {y_k}^2 
}
{2 \sigma_2^2}
] d {}{x_1} d {}{x_2}\cdots dx_nd {}{y_1} d {}{y_2}\cdots dy_m
\nonumber
\\
=
&
\frac{1}{{\sqrt{2 \pi }(\frac{\sigma_1^2}{n}+\frac{\sigma_2^2}{m})^{1/2}{}}}
\int_{{- \eta}}^{\eta} \exp[{}- \frac{{x}^2 }{2 (\frac{\sigma_1^2}{n}+\frac{\sigma_2^2}{m})}] d {x}
\label{eq62}
\end{align}
Solving the equation (\ref{eq28}),
we get that
\begin{align}
\eta^\gamma_{\omega} = 
(\frac{\sigma_1^2}{n}+\frac{\sigma_2^2}{m})^{1/2}
z(\frac{1-\gamma}{2})
\label{eq63}
\end{align}
Therefore,
for any 
$\widehat{x}$
$=$
$(x,y)$
$=(x_1,\ldots, x_n,y_1,\ldots, y_m)$
$({}\in {\mathbb R}^n\times {\mathbb R}^m)$,
we get $D_{\widehat{x}}^{\gamma}$
(
the $({}\gamma{})$-confidence interval
of
${\widehat x}$
)
as follows:
\begin{align}
D_{\widehat{x}}^{\gamma}
&
=
\{
{\omega}
(\in
\Omega)
:
d_\Theta ({}E(\widehat{x}),
\pi(\omega))
\le
\eta^\gamma_{\omega }
\}
\nonumber
\\
&
=
\{ (\mu_1, \mu_2 ) \in {\mathbb R} \times {\mathbb R}
\;:\;
|
\frac{\sum_{k=1}^n x_k }{n}
-
\frac{\sum_{k=1}^m y_k}{m}
-(\mu_1 - \mu_2 )|
\le 
(\frac{\sigma_1^2}{n}+\frac{\sigma_2^2}{m})^{1/2}
 z(\frac{1-\gamma}{2})
 \}
\label{eq64}
\end{align}
%
\subsection{The case that
$d^x_\Theta$ depends on $x$;
Student's t-distribution}
\par
\noindent
\rm
\par
The arguments in this section are continued from Example 2.

\par
\noindent
\bf
Example 7
\rm
[Student's t-distribution].
Consider the simultaneous measurement
${\mathsf M}_{C_0({\mathbb R} \times {\mathbb R}_+)}$
$({\mathsf O}_N^n = ({\mathbb R}^n, {\mathcal B}_{\mathbb R}^n, {{{N}}^n}) ,$
$S_{[(\mu, \sigma)]})$
in $C_0({\mathbb R} \times {\mathbb R}_+)$.
Thus,
we consider that
$\Omega = {\mathbb R} \times {\mathbb R}_+$,
$X={\mathbb R}^n$.
Put
$\Theta={\mathbb R}$ with the semi-distance
$d_\Theta^x (\forall x \in X)$
such that
\begin{align}
d_\Theta^x (\theta_1, \theta_2)
=
\frac{|\theta_1-\theta_2|}{{\overline{\sigma}'(x)}/\sqrt{n}}
\quad
\qquad
(\forall x \in X={\mathbb R}^n,
\forall \theta_1, \theta_2 \in \Theta={\mathbb R}
)
\label{eq65}
\end{align}
where ${\overline{\sigma}'(x)}=\sqrt{\frac{n}{n-1}}\overline{\sigma}(x)$.
The quantity $\pi:\Omega(={\mathbb R} \times {\mathbb R}_+)
\to
\Theta(={\mathbb R})$
is defined by
\begin{align}
\Omega(={\mathbb R} \times {\mathbb R}_+)
\ni \omega
=
(\mu, \sigma )
\mapsto \pi (\mu, \sigma )
=
\mu
\in
\Theta(={\mathbb R})
\label{eq66}
\end{align}
\rm
Also, define the estimator
$E:X(={\mathbb R}^n) \to \Theta(={\mathbb R})$
such that
\begin{align}
E(x)=E(x_1, x_2, \ldots , x_n )
=
\overline{\mu}(x)
=
\frac{x_1 + x_2 + \cdots + x_n}{n}
\label{eq67}
\end{align}
Let
$\gamma$
be a real number such that
$0 \ll \gamma < 1$,
for example,
$\gamma = 0.95$.
Thus, for any
$ \omega=(\mu, \sigma )  ({}\in\Omega=
{\mathbb  R} \times {\mathbb R}_+ )$,
we see that
\begin{align}
&
[N^n(\{ x \in X \;:\; 
d^x_\Theta ( E(x) , \pi( \omega ) )
\le \eta
\}
)](\omega )
\nonumber
\\
=&
[N^n(\{ x \in X \;:\; 
\frac{
|\overline{\mu}(x)- \mu |}{
{{\overline{\sigma}'(x)}/\sqrt{n}}
}
\le \eta
\}
)](\omega )
\nonumber
\\
=
&
\frac{1}{({{\sqrt{2 \pi }\sigma{}}})^n}
\underset{
- \eta
\le
\frac{
|\overline{\mu}(x)- \mu |}{
{{\overline{\sigma}'(x)}/\sqrt{n}}
}
\le \eta
}{\int \cdots \int}
\exp[{}- \frac{\sum_{k=1}^n ({}{}{x_k} - {}{\mu}  {})^2 
}
{2 \sigma^2}    {}] d {}{x_1} d {}{x_2}\cdots dx_n
\nonumber
\\
=
&
\frac{1}{({{\sqrt{2 \pi }{}}})^n}
\underset{
- \eta
\le
\frac{
|\overline{\mu}(x)- \mu |}{
{{\overline{\sigma}'(x)}/\sqrt{n}}
}
\le \eta
}
{\int \cdots \int}
\exp[{}- \frac{\sum_{k=1}^n ({}{}{x_k}  {}  {})^2 
}
{2 }    {}] d {}{x_1} d {}{x_2}\cdots dx_n
\nonumber
\\
\noindent
=
&
\int_{-\eta}^{\eta
}
p^t_{n-1}(x)
dx
\label{eq98}
\end{align}
where
$p^t_{n-1}$
is
the t-distribution with $n-1$ degrees of freedom.
Solving the equation
$
\gamma
=
\int_{-\eta^\gamma_{\omega}}^{\eta^\gamma_{\omega}
}
p^t_{n-1}(x)
dx
$,
we get $\eta^\gamma_{\omega}$
$=t((1-\gamma)/2)$.

%
%
%
%
%
%
%

Therefore,
for any
$x$
$({}\in X{})$,
we get $D_x^{\gamma}$(
\it
the $({}\gamma{})$-confidence interval
of
$x$
\rm
)
as follows:
\begin{align}
D_x^{\gamma}
&
=
\{
\pi({\omega})
(\in
\Theta)
:
\omega \in \Omega,
\;\;
d^x_\Theta ({}E(x),
\pi(\omega )
)
\le
\eta^\gamma_{\omega }
\}
\nonumber
\\
&
=
\{ \mu \in \Theta(={\mathbb R})
\;:\;
\overline{\mu}(x)
-
\frac{{\overline{\sigma}'(x)}}{\sqrt{n}}
t((1- \gamma)/2)
\le
\mu
\le
\overline{\mu}(x)
+
\frac{{\overline{\sigma}'(x)}}{\sqrt{n}}
t((1- \gamma)/2)
\}
\label{eq69} 
\end{align}

\rm

\section{Conclusions
}
\par
\noindent
\par
It is sure that
statistics and (classical) quantum language
are similar. however,
quantum language has the firm structure (\ref{eq1}),
i.e.,
\begin{align}
\underset{\mbox{(=MT(measurement theory))}}{\fbox{Quantum language}}
=
\underset{\mbox{(measurement)}}{\fbox{Axiom 1}}
+
\underset{\mbox{(causality)}}{\fbox{Axiom 2}}
+
\underset{\mbox{(how to use Axioms)}}{\fbox{linguistic interpretation}}
\label{eq70}
\end{align}
Hence,
as seen in this paper,
every argument cannot but become clear
in quantum language.
Thus, quantum language is suited for the theoretical arguments.
\par
\noindent
\par
In fact, Theorem 2 (the confidence interval methods n quantum language) says that
\begin{itemize}
\item{}
from the pure theoretical point of view, we can not understand the confidence interval methods
without the three concepts,
that is,
"estimator $E:X \to \Theta$" and "quantity ${\pi:\Omega \to \Theta}$" and "semi-distance $d^x_{\Theta}$",
\end{itemize}
which is also shown throughout Remarks 1-5
and Examples 4-7.
This is our new view-point of the confidence interval methods.
\vskip0.3cm
We hope that our approach will be examined from various points of view.


\rm
\par
\renewcommand{\refname}{
\large 
References}
{
\small

\normalsize
}

\end{document}